\documentclass{svmult}
\usepackage{graphicx, amsmath, amssymb}       
\usepackage[bottom]{footmisc}
\newcommand \bei {\begin{enumerate}}
\newcommand \eei {\end{enumerate}}
\newcommand \be         {\begin{equation}}
\newcommand \ee         {\end{equation}} 

\newcommand \bel {\be\label}
\newcommand \RR         {\mathbb{R}}
 
\newcommand \del        \partial
\newcommand \eps    \varepsilon 
\newcommand \Boxt {\widetilde \Box}
\newcommand \Rd {{R^\dag}}  
\newcommand \gd {{g^\dag}}  

\newcommand \Gammad { {\Gamma^{\dag}}}

\newcommand \Ocal   {\mathcal O}    
\newcommand \Hcal {\mathcal H}   
\newcommand \HH    {{\mathcal H}}

\begin{document}

\title*{The global nonlinear stability of Minkowski spacetime for self-gravitating massive fields. 
\\
A brief overview}
\titlerunning{The global nonlinear stability of Minkowski spacetime} 

\author{Philippe G. LeFloch\inst{1}}


\institute{Laboratoire Jacques-Louis Lions 
and 
Centre National de la Recherche Scientifique,  Universit\'e Pierre et Marie Curie, 4 Place Jussieu, 75258, Paris, France.
\\
Email: contact@philippelefloch.org}


\maketitle


\section{Introduction}

This is a short review of the series of papers \cite{PLF-YM-one,PLF-YM-two,PLF-YM-three} which, in collaboration with Yue Ma, establish several novel existence results for systems of coupled wave-Klein-Gordon equation. Our method --the Hyperbolic Hyperboloidal Method-- 
 has allowed us to address the global evolution problem for the Einstein equations of general relativity and investigate the global geometry of {\sl matter spacetimes} that are initially close to Minkowski spacetime. The Einstein equations  (when expressed in wave gauge) take the form of nonlinear system of partial differential equations of hyperbolic type and, in presence of self-gravitating massive matter fields, involve a strong coupling between wave equations (for the geometry) and Klein-Gordon equations (for the matter fields). Our method also provides a global existence theory for the field equations of the $f(R)$-theory of gravity, which is a natural generalization of Einstein's gravity theory (see below).

The global nonlinear stability problem for Minkowski spacetime is formulated from initial data which are prescribed on a spacelike hypersurface and
are a small perturbation of an asymptotically flat slice in Minkowski space. This problem is equivalent to a global-in-time existence problem for a system of nonlinear wave equations with sufficiently small data in a weighted Sobolev space. 

There are several major chalenges to be overcome. Gravitational waves are perturbations propagating in the curved spacetime and may relate to either the Weyl curvature (in  vacuum spacetimes) or the Ricci curvature (in matter spacetimes). It is required to understand the effect of nonlinear wave interactions on the possible growth of the energy, in order to be able to exclude dynamical instabilities and self-gravitating massive modes
and, therefore, to avoid gravitational collapse (trapped surfaces, black holes) --a generic phenomena in general relativity \cite{Christodoulou,Wald}. 
 
The global dynamics is particularly complex in general, but {\sl sufficiently small} perturbations in Minkowski spacetime 	are expected to disperse in timelike directions and an asymptotic convergence to Minkowski spacetime to be observed. More precisely, in order to prove its stability, it is necessary to establish that the spacetime is future timelike geodesically complete. 

We will begin by reviewing Einstein's gravity and the f(R)-gravity theory and explain their relation. Next, the wave-Klein-Gordon formalism will be presented, and finally the global nonlinear stability will be stated. For further reading, we refer to the works by
Donninger and Zenginoglu \cite{DZ},  
Fajman, Joudioux, and Smulevici \cite{FJS}, Wang \cite{Wang}, and Zenginoglu \cite{Zenginoglu}.

  
\section{Self-gravitating massive fields} 

Throughout, we are interested in Lorentzian manifolds (or spacetime) $(M, g_{\alpha\beta})$ with signature $(-, +, +, +)$ and in local coordinates we write $g= g_{\alpha\beta} dx^\alpha dx^\beta$. 
For instance, Minkowski spacetime $M=\RR^{3+1}$ is described in standard coordinates by $g_M = -(dx^0)^2 + \sum_{a=1}^3 (dx^a)^2$.  
In a spacetime, the covariant derivative operator allows to write schematically $\nabla_\alpha X = \del_\alpha X + \Gamma \star X$ with $\Gamma \simeq \del g$. (The exact expressions in coordinates will be given only later in this text.) 
The Ricci curvature also reads schematically $R_{\alpha\beta} = \del^2 g + \del g \star \del g$ and   
one also defines the scalar curvature $R := R_\alpha^\alpha = g^{\alpha\beta} R_{\alpha\beta}$ by taking the trace of the Ricci curvature. 
Here, $\alpha, \beta =0, 1, 2, 3$ and, whener relevant, Einstein's summation convention on repeated indices is in order. 

The Einstein equations for self-gravitating matter have the form 
\be
G_{\alpha\beta} = 8\pi T_{\alpha \beta},
\ee
in which $G_{\alpha\beta} := R_{\alpha\beta} - (R/2) g_{\alpha\beta}$ is the Einstein curvature tensor and $T_{\alpha \beta}$ denotes the energy-momentum of the matter.
A  (minimally coupled) massive scalar field with potential $U(\phi)$, for instance with the quadratic potential 
$U(\phi) = \frac{c^2}{2} \phi^2$, is described by the energy-momentum tensor
\be
T_{\alpha\beta} := \nabla_\alpha \phi \nabla_\beta \phi - \Big( {1 \over 2}
g^{\alpha'\beta'} \nabla_{\alpha'} \phi \nabla_{\beta'} \phi + U(\phi) \Big) g_{\alpha\beta}.
\ee
Consequently, the Einstein-Klein-Gordon system for the unknown $(M, g_{\alpha\beta}, \phi)$ reads 
\be
\aligned
R_{\alpha\beta} - 8\pi \Big(\nabla_\alpha\phi\nabla_{\beta} \phi + U(\phi) \, g_{\alpha\beta} \Big) & = 0,
\\
\Box_g\phi - U'(\phi) & = 0,
\endaligned
\ee
where $\Box_g =\nabla_\alpha\nabla^\alpha$ denotes the wave operator associated with the unknown metric $g$.  The above equation is a geometric PDE's, which enjoys a gauge invariance property.  
 

On the other hand, the field equations for the $f(R)$-theory of modified gravity are based on the following generalized Hilbert-Einstein functional: 
\be
\int_M \Big( f(R) + 16 \pi L[\phi, g]\Big) \, dV_g, 
\ee
in which the nonlinear function  $f(R) = R + \frac{\kappa}{2} R^2 + \kappa^2 \Ocal( R^3)$ is prescribed with $\kappa>0$. The condition $\kappa:=f''(0) >0$ will be essential for global stability. This theory has a long history in physics, beginning with Weyl (1918), Pauli (1919), Eddington (1924), and many others 
\cite{BransDicke, Buchdahl}.  
We emphasize that alternative theories of gravity are relevant in view of recent observational data,
which have demonstrated the accelerated expansion of the Universe and have identified instabilities in galaxies in our Universe. The $f(R)$-theory allows for the gravitation field to be mediated by an additional field without explicitly introducing a notion of 'dark matter'. 
  
Numerical evidence and physical heuristics have led to the conjecture that asymptotically flat, matter spacetimes should be stable \cite{OCP}, even though the existence of a family of ``oscillating soliton stars'' had first suggested a possible instability mechanism within small perturbations of massive fields. Advanced numerical methods were necessary to handle the long-time evolution of oscillating soliton stars. During an initial phase, the matter {\sl tends to collapse}, but during an intermediate phase (below a certain threshold in the mass density) the {\sl collapse slows down}, until finally the {\sl dispersion} becomes of the main feature of the evolution of the matter field.  

Observe that in asymptotically anti-Sitter (AdS) spacetimes, such instabilities are observed and the effect of gravity is dominant so that generic (even arbitrarily small) initial data lead to black hole formation. In AdS spacetime, the matter is confined and cannot disperse: the timelike boundary is reached in finite proper time

  
\section{The wave-Klein-Gordon formulation}

Unless we introduce a specific gauge, the field equations $G_{\alpha\beta} = 8\pi T_{\alpha\beta}$ in coordinates take the form of a second-order system with no specific PDE type. By imposing the wave gauge $\Box_g x^\gamma = 0$, we find 
\be
2 \, 
g^{\alpha\beta} \, \del_\beta g_{\alpha\gamma} - g^{\alpha\beta} \, \del_\gamma g_{\alpha\beta} = 0, 
\qquad \gamma = 0, \ldots, 3, 
\ee
leading us to an expression of the Ricci curvature $R_{\alpha\beta} \simeq \Box_g g_{\alpha\beta}$ 
(after Einstein, Choquet-Bruhat, De Turck, etc.). The Einstein-massive field system is then equivalent to a second-order system of $11$ nonlinear wave-Klein-Gordon equations, supplemented with the Hamiltonian-momentum Einstein's contraints.  
Namely, in wave gauge, the Einstein equations for a self-gravitating massive field read 
\bel{eq:120} 
\aligned 
\widetilde {\Box}_g g_{\alpha\beta}
= \, &  F_{\alpha\beta}(g, \del g)  
- 8 \pi \, \big( 2  \del_\alpha \phi\del_\beta \phi + c^2 \phi^2 \, g_{\alpha \beta} \big), 
\\
\widetilde {\Box}_g \phi  - c^2 \phi 
= \, &  0,
\endaligned
\ee
with $\Boxt_g \psi:= g^{\alpha' \beta'}\del_{\alpha'}\del_{\beta'} \psi$.

The expression of the quadratic nonlinearities $F_{\alpha\beta}(g, \del g)$ is given in the following lemma and involves null terms of the general form $g^{\alpha\beta} \del_\alpha u \del_\beta u$ or $\del_\alpha u \del_\beta v-  \del_\beta u \del_\alpha v$, as well as terms that do not satisfy the null ondition and requires a specific analysis. Recall first that 
$$
\aligned
R_{\alpha\beta} & = \del_{\lambda} \Gamma_{\alpha\beta}^{\lambda} - \del_{\alpha} \Gamma_{\beta\lambda}^{\lambda} + \Gamma_{\alpha\beta}^{\lambda} \Gamma_{\lambda\delta}^{\delta} - \Gamma_{\alpha\delta}^{\lambda} \Gamma_{\beta\lambda}^{\delta}
\\  
\Gamma_{\alpha\beta}^{\lambda} 
& = {1 \over 2} \gd^{\lambda\lambda'} \big(\del_{\alpha}\gd_{\beta\lambda'}
+ \del_{\beta}\gd_{\alpha\lambda'} - \del_{\lambda'}\gd_{\alpha\beta} \big)
\endaligned
$$   
Following Lindblad and Rodnianski \cite{LR2} for vacuum Einstein spacetime, we have the following result.

\begin{lemma} With $F_{\alpha\beta} = Q_{\alpha\beta} + P_{\alpha\beta}$, the Ricci curvature in wave gauge reads 
$$
\aligned
2 \, R_{\alpha\beta} 
& = - \Boxt_g g_{\alpha\beta} 
  + Q_{\alpha\beta} + P_{\alpha\beta},  
\endaligned
$$ 
which contains  
\bei

\item null terms satisfying Klainerman's null condition (and enjoying good decay in time)
\be
\aligned
Q_{\alpha \beta}:& = 
  g^{\lambda \lambda'} g^{\delta \delta'} \del_{\delta} g_{\alpha \lambda'} \del_{\delta'} g_{\beta \lambda}
\\
&     \quad - g^{\lambda \lambda'} g^{\delta \delta'} \big
(\del_{\delta} g_{\alpha \lambda'} \del_{\lambda} g_{\beta \delta'} - \del_{\delta} g_{\beta \delta'} \del_{\lambda} g_{\alpha \lambda'} \big)
\\
& \quad +  g^{\lambda \lambda'} g^{\delta \delta'}
\big(\del_\alpha g_{\lambda'\delta'} \del_{\delta} g_{\lambda \beta} - \del_\alpha g_{\lambda \beta} \del_{\delta} g_{\lambda'\delta'} \big) + \ldots,  
\endaligned
\ee

\item and quasi-null terms (as they are called by the authors)
\be
P_{\alpha \beta} := -  {1 \over 2}  g^{\lambda \lambda'} g^{\delta \delta'} \del_\alpha g_{\delta \lambda'} \del_\beta g_{\lambda \delta'}
+ {1 \over 4} g^{\delta \delta'} g^{\lambda \lambda'} \del_\beta g_{\delta \delta'} \del_\alpha g_{\lambda \lambda'} 
\ee
(which will require a further investigation based on the wave gauge condition). 
\eei  
\end{lemma}

A similar decomposition can be written for the conformal metric of the f(R)-theory of gravity, which we now introduced. 
The modified gravity equations  read 
\be
N_{\alpha\beta} = 8\pi T_{\alpha\beta}
\ee
and 
are based on a choice of a function $f(R) = R + {\kappa \over 2} R^2 + \ldots$, and take the form of 
a fourth-order system with no specific PDE type. We propose to rely on an augmented formulation with unknown $(\gd_{\alpha\beta}, \rho)$ defined as follows, by regarding the spacetime curvature as an independent unknown and by working with the conformal metric
\be
\gd_{\alpha\beta}:= f'(R_g) g_{\alpha\beta}, 
\ee
in which $\rho := {1 \over \kappa} \ln f'(R_g)$. In view of the standard 
relation between the Ricci curvature tensors of $g$ and $\gd$, i.e. 
$$
\Rd_{\alpha\beta} = R_{\alpha\beta} - 2\big(\nabla_\alpha\nabla_{\beta}\rho
- \nabla_\alpha\rho\nabla_{\beta}\rho\big) - \big(\Box_g\rho+2g(\nabla \rho,\nabla \rho)\big)g_{\alpha\beta}, 
$$ 
we arrive at a third-order system. In addition, from the trace of the field equation, we derive an evolution equation for the scalar curvature 
which is a new degree of freedom in the theory and must be supplemented with suitable initial data. Finally, in wave coordinates
\be
\Box_\gd x^\alpha = 0
\ee 
we arrive at a second-order system of $12$ nonlinear wave-Klein-Gordon equations. The structure is analogous to the one of the 
Einstein-massive field system, but has a significantly more involved algebraic structure and admits additional constraints.

\begin{proposition} The equations of f(R)-gravity for a self-gravitating massive field take the form 
\be
\aligned 
\widetilde {\Box}_\gd  \gd_{\alpha\beta}
= \, &  F_{\alpha\beta}(\gd,\del \gd) 
- 8 \pi \, \big( 2 e^{-\kappa \rho} \del_\alpha \phi\del_\beta \phi + c^2 \phi^2 e^{-2\kappa \rho}  \, \gd_{\alpha \beta} \big)
\\  
&  \hskip2.cm
 - 3 \kappa^2 \del_{\alpha}\rho\del_{\beta}\rho+  \kappa \, {\mathcal O}(\rho^2)
\gd_{\alpha\beta}, 
\\ 
\widetilde {\Box}_\gd \phi  - c^2 \phi 
= \,  
&  c^2 \big( e^{-\kappa\rho} - 1 \big) \phi   + \kappa \gd^{\alpha\beta}\del_{\alpha}\phi\del_{\beta}\rho,
\\
3 \kappa \, \widetilde {\Box}_\gd \rho - \rho 
= \,  
& \kappa \, {\mathcal O}(\rho^2)
- 8 \pi\Big( \gd^{\alpha\beta}  \del_\alpha \phi \del_\beta \phi + {c^2 \over 2} \, e^{-\kappa\rho} \phi^2 \Big),
\endaligned
\ee
supplemented with: 
\bei 

\item the wave gauge conditions $\gd^{\alpha\beta} \Gammad_{\alpha\beta}^{\lambda} = 0$, 

\item the curvature compatibility condition  $e^{\kappa\rho} = f'(R_{e^{-\kappa\rho} \gd})$,   

\item and the Hamiltonian and momentum constraints. 

\eei
These three sets of conditions can be propagated from a Cauchy hypersurface.
\end{proposition}

\begin{proposition} In the limit $\kappa \to 0$ one finds
\be
\gd \to g
\quad \text{ and } 
\rho \to  8 \pi\big( g^{\alpha\beta} \nabla_\alpha \phi \nabla_\beta \phi + {c^2 \over 2} \phi^2 \big), 
\ee
and the Einstein system for a self-gravitating massive field is \eqref{eq:120}.  
\end{proposition}
   
This completes the formulation of the field equations in a PDE form and the geometric problem of interest can be reformulated as a global existence problem for coupled nonlinear wave equations. 
Our main challenge is that the system is {\sl not invariant by scaling},  
and one must rely on {\sl fewer symmetries} 
in, for instance, defining weighted energy-like functionals. 
The analysis of coupled wave equations and Klein-Gordon equations is particularly challenging
and drastically {\sl different time asymptotic behavior} arise for the unknown components of the system: $O(t^{-1})$ for wave equations and $O(t^{-3/2})$ for Klein-Gordon equations. 

We also need to investigate the dependence in $f$ and determined the {\sl singular limit} $f(R) \to R$ which, as will can see, transforms a second-order PDE into an algebraic equation.
 
  
\section{The global nonlinear stability} 

For the global existence theory, we need to establish that there is a sufficient rate of time decay for all the noninearities of interest. 
The nonlinear coupling between the geometry and massive matter 
leads to strong interactions at the PDE level and, consequently, it is necessary to 
be able to establish (almost) sharp $L^2$ and $L^\infty$ time-decay for the metric and matter field. 
Understanding the quasi-null structure of the Einstein equations is fundamental since the standard 
null condition is violated and an amplification phenomena arise for the energy. 
 
We thus consider the initial value problem for the Einstein equations (and its generalization). An initial data set, by definition,  provides us with the geometry of the initial hypersurface $(M_0 \simeq \RR^3, g_0, k_0)$
and the initial data for the matter field $\phi_0, \phi_1$. We assume that these data  
are sufficiently close to a spacelike, asymptotically flat slice in Minkowski spacetime. The local existence is standard and goes back to  
Choquet-Bruhat \cite{YCB}: to each initial data set, one can assciate a unique maximal, globally hyperbolic Cauchy development
(i.e., intuitively, the maximal part of the spacetime which is uniquely determined by the prescribed initial data and remains smooth). 
 
The fundamental work on the stability of Minkowski spacetime for the vacuum Einstein equations (or massless matter)
was done by Christodoulou and Klainerman \cite{CK} (and later generalized in \cite{BieriZipser}):  

\bei 

\item They introduced  a fully geometric proof, in which the Bianchi identities are regarded as the main evolution equations, 

\item they analyzed the geometry of null cones and defined a double null-maximal foliation,

\item and they relied on all of the Killing fields of Minkowski spacetime.

\eei 

The earlier work by Friedrich \cite{Friedrich83} also addressed the global existence problem for the vacuum Einstein equations, and established the nonlinear stability of DeSitter spacetime. More recently, Lindblad and Rodnianski \cite{LR2} 
obtained the first global existence result for the vacuum Einstein equations in wave coordinates (despite an ``instability'' result by Choquet-Bruhat) and again relied on all of the Killing fields of Minkowski spacetime. 
They introduced a foliation by asymptotically flat hypersurfaces. 
 
In constrast, the recent work  \cite{PLF-YM-one,PLF-YM-two,PLF-YM-three}  addresses this stability problem for {\sl self-gravitating massive matter fileds}:  

\bei 
 
\item The proposed new method (the Hyperboloidal Foliation Method) 
does not rely on Minkowski's scaling field $r \del_r + t \del_t$, 

\item which is the key of be able to tackle massive matter fields, 

\item is based on an asymptotically hyperbolic foliation,  

\item and leads to a somewhat simpler proof for the case of massless fields. 

\eei 
 
The positive mass theorem restricts the possible behavior of solutions at spacelike infinity. No solution can be exactly Minkowski ``at infinity'', but can coincide with the Schwarzschild metric outside a spatially compact region. More generally, solutions are assumed to 
approach the Schwarzschild metric near space infinity (with ADM mass $m <<1$). 
We only provide here informal statements of our results and we refer to  \cite{PLF-YM-one,PLF-YM-two,PLF-YM-three}  for the precise statements. 

\begin{theorem}[Nonlinear stability of Minkowski spacetime with self-gravitating massive fields]
Consider the Einstein-massive field system when the initial data set $(M_0 \simeq \RR^3, g_0, k_0, \phi_0, \phi_1)$ is 
asymptotically Schwarz\-schild and sufficiently close to Minkowski data and satisfies the Einstein constraint equations. 
Then, the initial value problem 

\bei 

\item admits a globally hyperbolic Cauchy development, 

\item which is foliated by asymptotically hyperbolic hypersurfaces, 

\item and is future causally geodesically complete and asymptotically approaches Minkowski spacetime.  

\eei 
\end{theorem}
 
\begin{theorem}[Nonlinear stability of Minkowski spacetime in  f(R)-gravity]
Consider the field equations of $f(R)$-modified gravity when the initial data set 
$(M_0 \simeq \RR^3, g_0, k_0, R_0, R_1, \phi_0, \phi_1)$ is 
 asymptotically Schwarzschild and sufficiently close to Minkowski data
and satisfies the constraint equations of modified gravity.  
Then, the initial value problem 

\bei 

\item admits a globally hyperbolic Cauchy development,

\item  which is foliated by asymptotically hyperbolic hypersurfaces,

\item and is future causally geodesically complete and asymptotically approaches Minkowski spacetime. 

\eei 
\end{theorem} 

The limit problem $\kappa \to 0$ can be viewed as a relaxation phenomena for the spacetime scalar curvature. We pass from the second-order wave equation
\be
3 \kappa \, \widetilde {\Box}_\gd \rho - \rho 
=  \kappa \, {\mathcal O}(\rho^2) 
- 8 \pi\Big( \gd^{\alpha\beta}  \del_\alpha \phi \del_\beta \phi + {c^2 \over 2} \, e^{-\kappa\rho} \phi^2 \Big) 
\ee
 to the purely algebraic equation
\be
\rho \to  8 \pi\big( g^{\alpha\beta} \nabla_\alpha \phi \nabla_\beta \phi + {c^2 \over 2} \phi^2 \big). 
\ee
 
\begin{theorem}[f(R)-spacetimes converge toward Einstein spacetimes]
In the limit $\kappa \to 0$, when the nonlinear function $f=f(R)$ (the integrand in the Hilbert-Einstein action)
 approaches the scalar curvature function $R$, the Cauchy developments of modified gravity (given in the previous theorem) 
converge (in every bounded time interval, in a sense specified quantitatively in Sobolev norms) to Cauchy developments of Einstein's gravity theory.
\end{theorem}

The proofs rely on weighted norms associated with the asymptotically hyperboloidal foliation which we construct. 
Our energy norms are solely based on the translations $\del_\alpha$ and the Lorentzian boosts
$L_a$ of Minkowski spacetime. These fields enjoy good commutator properties even in curved space
 and allow us to decompose the wave operators, the metric, etc.  
On each hyperbooidal hypersurface $\HH^n[s]$ at any hyperboloidal time $s$, in wave coordinates, we use the boosts to define  the norm 
\be
\aligned
\big( 
\| u \|_{\HH^n[s]}\big)^2  
: =  
\sup_{a=1, 2, 3}
\sum_{|J|\leq n}  \int_{\Hcal_s \simeq \RR^3} |L_a^J  u |^2 \, dx 
\endaligned
\ee
and, within the spacetime, we use the translations to define the norm 
\be
\|u\|_{\HH^N[s_0, s_1]} := \displaystyle
\sup_{s \in [s_0, s_1]} \sum_{|I| +n \leq N}  
\big\|\del^I u \big\|_{\HH^n[s]}. 
\ee
We introduce a suitable bootstrap argument, which shows that the total contribution of the interaction terms contributes only a finite amount to the growth of the total energy. We derive global time-integrability properties for the source terms, which are established from 
sharp pointwise estimates ---required to handle the strong geometry-matter interactions under consideration. 
Sobolev inequalities and Hardy inequalities are adapted to the hyperboloidal foliation, and a hierarchy of energy bounds distinguishes between various orders of differentiation and growth rates in the hyperboloidal time $s$.
 

\end{document}